\documentclass[12pt]{amsart}
\theoremstyle{plain}
\newtheorem{theorem}{Theorem}[section]
\newtheorem{corollary}[theorem]{Corollary}

\newtheorem{proposition}[theorem]{Proposition}

\theoremstyle{definition}

\theoremstyle{remark}

\newtheorem{problem}[theorem]{Problem}

\begin{document}

\title[C*-algebras of compact operators and Hilbert C*-modules]{Characterizing
C*-algebras of compact operators by generic categorical
properties of Hilbert C*-modules}
\author[M.~Frank]{Michael Frank}
\dedicatory{{Dedicated to the memory of Yu.~P.~Solovyov}}
\address{Hochschule f\"ur Technik, Wirtschaft und Kultur (HTWK) Leipzig,
Fachbereich IMN, Gustav-Freytag-Strasse 42A, D-04277 Leipzig,
 \linebreak[4] Germany}
\email{mfrank@imn.htwk-leipzig.de}
\begin{abstract}
B.~Magajna and J.~Schweizer showed in 1997 and 1999, respectively, that
C*-algebras of compact operators can be characterized by the property that
every norm-closed (and coinciding with its biorthogonal complement, resp.)
submodule of every Hilbert C*-module over them is automatically
an orthogonal summand. We find out further generic properties of the category
of Hilbert C*-modules over C*-algebras which characterize precisely the
C*-algebras of compact operators.
\end{abstract}
\maketitle

In 1997 B.~Magajna obtained the equivalence of the property of the
category of Hilbert C*-modules over a certain C*-algebra $A$ that
any Hilbert C*-submodule is automatically an orthogonal summand
with the property of the C*-algebra $A$ of coefficients to admit a
faithful $*$-representation in some C*-algebra of compact
operators on some Hilbert space, cf.~Theorem \ref{Mag-Schw}. In
1999 J.~Schweizer was able to sharpen the argument replacing the
Hilbert C*-module property of B.~Magajna by the property of the
category of Hilbert C*-modules over a certain C*-algebra $A$ that
any Hilbert C*-submodule which coincides with its biorthogonal
complement is automatically an orthogonal summand, cf.~Theorem
\ref{Mag-Schw}. Later on in 2003 M.~Kusuda published further
results which indicate that in the majority of situations the
Hilbert C*-module property can be weakened merely requiring the
${\rm K}_A({\mathcal M})$-$A$-subbimodules of the Hilbert
C*-modules $\mathcal M$ to be always orthogonal summands, see
\cite{Kusuda,Kusuda2} for the details.

Studying the work of B.~Magajna and J.~Schweizer
C*-algebras $A$ of the form $A= c_0$-$\sum_\alpha \oplus {\rm K}(H_\alpha)$
become of special interest, where the symbol ${\rm K}(H_\alpha)$ denotes the
C*-algebra of all compact operators on some Hilbert space $H_i$, and the
$c_0$-sum is either a finite block-diagonal sum or a block-diagonal sum with a
$c_0$-convergence condition on the C*-algebra components ${\rm K}(H_\alpha)$.
The $c_0$-sum may possess arbitrary cardinality. This kind of C*-algebras has
been precisely characterized by W.~Arveson \cite[\S I.4, Th.~I.4.5]{Arveson}
as the C*-subalgebras of (full) C*-algebras of compact operators
on Hilbert spaces. Throughout the present paper we refer to these C*-algebras
as to {\it C*-algebras of compact operators on certain Hilbert spaces}.

In 1998-99 V.~I.~Paulsen and the author investigated injective and projective
objects in categories of Hilbert C*-modules over certain C*-algebras,
cf.~\cite{FP}. The research work revealed a number of further properties
every Hilbert C*-submodule of any Hilbert C*-module over a fixed C*-algebras
$A$ might admit. The goal of the present paper is to present such results
which all lead to further equivalent characterizations of the C*-algebras of
compact operators. Together with results of D.~Baki\'c and B.~Gulja{\v{s}}
(\cite{BG}) on the existence of modular normalized tight frames in Hilbert
C*-modules over C*-algebras of compact operators we collect a large number
of useful properties of these classes of Hilbert C*-modules.
(For a concise introduction to modular frame theory we refer to \cite{FL}.)

The paper consists of a section with introductory material and of a second
section which contains the results.

\section{Preliminaries}

We give definitions and basic facts of Hilbert C*-module
theory needed for our investigations. The papers
\cite{Pa1,Kas,Fr90,Lin:90/2,Lin1,Fr98}, some chapters in \cite{JT,NEWO}, and
the books by E.~C.~Lance \cite{Lance} and by I.~Raeburn and D.~P.~Williams
\cite{RaeWil} are used as standard sources of reference. We make the
convention that all C*-modules of the present paper are left modules by
definition. A {\it pre-Hilbert $A$-module over a C*-algebra} $A$ is an
$A$-module $\mathcal M$ equipped with an $A$-valued mapping $\langle .,.
\rangle : {\mathcal M} \times {\mathcal M} \rightarrow A$ which is $A$-linear
in the first argument and has the properties:
\[
\langle x,y \rangle = \langle y,x \rangle^* \; \, , \: \;
\langle x,x \rangle \geq 0 \quad {\rm with} \: {\rm equality} \: {\rm iff}
\quad x=0 \, .
\]
The mapping $\langle .,. \rangle$ is said to be {\it the $A$-valued inner
product on} $\mathcal M$. A pre-Hilbert $A$-module $\{ \mathcal M, \langle
.,. \rangle \}$ is {\it Hilbert} if and only if it is complete with respect
to the norm $\| . \| = \| \langle .,. \rangle \|^{1/2}_A$. We always assume
that the linear structures of $A$ and $\mathcal M$ are compatible. Two Hilbert
$A$-modules are {\it isomorphic} if they are isometrically isomorphic as
Banach $A$-modules, if and only if they are unitarily isomorphic,
\cite{Lance}. We would like to point out that Banach $A$-modules can carry
unitarily non-isomorphic $A$-valued inner products which induce equivalent
norms to the given one, nevertheless, \cite{Fr98}.

Hilbert C*-submodules of Hilbert C*-modules might not be direct summands, and
if they are direct summands then they can be topological, but not orthogonal
summands. We say that a Hilbert C*-module $\mathcal N$ is a topological
summand of a Hilbert C*-module $\mathcal M$ which contains $\mathcal N$ as a
Banach C*-submodule in case $\mathcal M$ can be decomposed into the direct sum
of the Banach C*-submodule $\mathcal N$ and of another Banach C*-submodule
$\mathcal K$. The denotation is ${\mathcal M} = {\mathcal N} \stackrel{.}{+}
{\mathcal K}$. If, moreover, the decomposition can be arranged as an orthogonal
one (i.e.~$\mathcal N \bot \mathcal K$) then the Hilbert C*-submodule
$\mathcal N \subseteq \mathcal M$ is an orthogonal summand of the Hilbert
C*-module $\mathcal M$. Examples of any kind of appearing situations can be
found in \cite{Fr98}.

Finally, we are going to consider various bounded C*-linear operators $T$
between Hilbert C*-modules $\mathcal M$, $\mathcal N$ with one and the same
C*-algebra of coefficients. Quite regularly those operators $T$ may not
admit an adjoint bounded C*-linear operator $T^*: \mathcal N \to \mathcal M$
fulfilling the equality $\langle T(x),y \rangle_{\mathcal N} = \langle x,
T^*(y) \rangle_{\mathcal M}$ for any $x \in \mathcal M$, any $y \in \mathcal
N$. We denote the C*-algebra of all bounded C*-linear adjointable operators
on a given Hilbert $A$-module $\mathcal M$ by ${\rm End}_A^*(\mathcal M)$.
The Banach algebra of all bounded $A$-linear operators on $\mathcal M$ is
denoted by ${\rm End}_A(\mathcal M)$. For more detailed information on
such situations we refer to \cite{Fr98}.

\section{C*-algebras of compact operators and the Magajna-Schweizer theorem}

This section aims to characterize the class of C*-algebras of
compact operators on certain Hilbert spaces and their
C*-subalgebras by the appearance of certain properties common to
all Hilbert C*-modules over them. Our starting point is the
following result by B. Magajna and J. Schweizer and some of its
immediate consequences:

\begin{theorem} \label{Mag-Schw} {\rm (B.~Magajna, J.~Schweizer
                     \cite{Mag,Schweiz})} \newline
   Let $A$ be a C*-algebra. The following three conditions are equivalent:

   \newcounter{cou001}
   \begin{list}{(\roman{cou001})}{\usecounter{cou001}}
   \item  $A$ is of $c_0$-$\sum_i \oplus {\rm K}(H_i)$-type, i.e.~it has a
          faithful $*$-representation as a  C*-algebra of compact operators
          on some Hilbert space.
   \item  For every Hilbert $A$-module $\mathcal M$ every Hilbert $A$-submodule
          \linebreak[4]
          $\mathcal N \subseteq \mathcal M$ is automatically orthogonally
          complemented in $\mathcal M$, i.e.~$\mathcal N$ is an orthogonal
          summand of $\mathcal M$.
   \item  For every Hilbert $A$-module $\mathcal M$ every Hilbert $A$-submodule
          \linebreak[4]
          $\mathcal N \subseteq \mathcal M$ that coincides with its
          bi-orthogonal complement ${\mathcal N}^{\bot\bot} \subseteq \mathcal
          M$ is automatically orthogonally complemented in $\mathcal M$.
   \end{list}
\end{theorem}

\begin{corollary} \label{Fr-neu}
   Let $A$ be a C*-algebra. The following three conditions are equivalent:

   \begin{list}{(\roman{cou001})}{\usecounter{cou001}}
   \item  $A$ is of $c_0$-$\sum_i \oplus {\rm K}(H_i)$-type, i.e.~it has a
          faithful $*$-representation as a  C*-algebra of compact operators
          on some Hilbert space.
   \setcounter{cou001}{3}
   \item For every Hilbert $A$-module $\mathcal M$ and every bounded $A$-linear
         map $T: {\mathcal M} \to {\mathcal M}$ there exists an adjoint bounded
         $A$-linear map $T^*:{\mathcal M} \to {\mathcal M}$.
   \item For every pair of Hilbert $A$-modules $\mathcal M$, $\mathcal N$ and
         every bounded $A$-linear map $T: {\mathcal M} \to {\mathcal N}$ there
         exists an adjoint bounded $A$-linear map $T^*:{\mathcal N} \to
         {\mathcal M}$.
   \end{list}
\end{corollary}

\begin{proof}
Condition (i) is equivalent to the assertion that every Hilbert $A$-module
$\mathcal M$ is orthogonally complemented as a Hilbert $A$-submodule of
arbitrary Hilbert $A$-modules by Theorem \ref{Mag-Schw}. The latter condition
on a certain $\mathcal M$ is equivalent to the assertion (iv) of the corollary
by \cite[Th.~6.3]{Fr98} - in case the range of the C*-valued inner product
is dense in the C*-algebra of coefficients. However, the range of the C*-valued
inner product of Hilbert C*-modules over C*-algebras of type (i) gives always
rise to an ideal of them, and these norm-closed two-sided ideals are of type (i)
again. So the restriction does not matter, and the equivalence of (i) and (iv)
follows.

To establish the equivalence of (i) with the last condition consider the
operator $T$ described at (v) as a bounded $A$-linear operator on the Hilbert
$A$-module ${\mathcal M} \oplus {\mathcal N}$ mapping pairs $(x,y)$ to pairs
$(0,T(x))$. Then the equivalence (i)$\leftrightarrow$(iv) applied to this
particular situation implies condition (v) as a simple calculation shows. The
converse conclusion becomes trivial resorting to the situation ${\mathcal M}
= {\mathcal N}$.
\end{proof}

The next group of facts is concerned with particular Hilbert C*-modules,
bounded module operators on them and their properties. We aim to describe
C*-algebras $A$ of compact operators on Hilbert spaces by
general properties of kernels of bounded module operators on Hilbert
$A$-modules and of images of bounded module operators on Hilbert C*-modules
that admit a closed range. By the way we prove the closed graph theorem for
bounded module operators on Hilbert C*-modules.

\begin{proposition} \label{closed-range} {\rm (N.~E.~Wegge-Olsen
                          \cite[Th.~15.3.8]{NEWO})} \newline
   Let $A$ be a C*-algebra, $\{ {\mathcal M}, \langle .,. \rangle \}$ be a Hilbert
   $A$-module and $T$ be an adjointable bounded module operator on $\mathcal M$.
   If $T$ has closed range then $T^*$, $(T^*T)^{1/2}$ and $(TT^*)^{1/2}$ have
   also closed ranges and
   \begin{eqnarray*}
   {\mathcal M} & = & {\rm Ker}(T) \oplus T^*({\mathcal M})
              =   {\rm Ker}(T^*) \oplus T({\mathcal M}) \\
            & = & {\rm Ker}(|T|) \oplus |T|({\mathcal M})
              =   {\rm Ker}(|T^*|) \oplus |T^*|({\mathcal M}) \, .
   \end{eqnarray*}
   In particular, each orthogonal summand appearing on the right is
   automatically norm-closed and coincides with its bi-orthogonal complement
   inside $\mathcal M$. Moreover, $T$ and $T^*$ have polar decomposition.
\end{proposition}

\begin{corollary}  {\rm (bounded closed graph theorem)} \label{adjointability}
    \newline
   Let $A$ be a C*-algebra and $\{ \mathcal M, \langle .,. \rangle \}$,
   $\{ \mathcal N, \langle .,. \rangle \}$ be two Hilbert $A$-modules. A
   bounded $A$-linear operator $T: \mathcal M \to \mathcal N$ possesses an
   adjoint operator $T^*: \mathcal N \to \mathcal M$ if and only if the
   graph of $T$ is an orthogonal summand of the Hilbert $A$-module $\mathcal M
   \oplus {\mathcal N}$. Beside this equivalence, the graph of every bounded
   $A$-linear operator $T$ coincides with its bi-orthogonal complement in
   $\mathcal M \oplus \mathcal N$, and it is always a topological summand
   with topological complement $\{ (0,z) : z \in \mathcal N \}$.
   By a counterexample due to E.~C.~Lance (\cite[pp.~102-104]{Lance}) this
   fails for some closed, self-adjoint, densely defined, unbounded module
   operators on certain Hilbert C*-modules.
\end{corollary}

\begin{proof}
Since the inequality $\|T(x)\| \leq \|T\| \|x\|$ is valid for every $x \in
\mathcal M$ the graph of $T$ is a norm-closed Hilbert $A$-submodule of the
Hilbert $A$-module $\mathcal M \oplus \mathcal N$. Moreover, since the graph
of $T$ is the kernel of the bounded module operator $S: (x,y) \to (0,T(x)-y)$
on $\mathcal M \oplus \mathcal N$ it coincides with its bi-orthogonal
complement there, \cite[Cor.~2.7.2]{Fr97}. If $T$ has an adjoint then the
operator $T':(x,y) \to (x,T(x))$ is adjointable on $\mathcal M \oplus
\mathcal N$. By Proposition \ref{closed-range} the graph of $T$ is an
orthogonal summand.

Conversely, if the graph of $T$ is an orthogonal summand of ${\mathcal M}
\oplus {\mathcal N}$ then its orthogonal complement consists precisely
of the pairs of elements $\{ (x,y) : x = -T^*(y), y \in {\mathcal N} \}$
since $\langle z,x \rangle_{\mathcal M} + \langle T(z),y \rangle_{\mathcal N}=0$
forces $T^*(y)(z) = \langle z,(-x) \rangle_{\mathcal M}$ for any $z \in
\mathcal M$ and $T^*:{\mathcal N} \to {\mathcal M}'$. So $T^*$ is everywhere
defined on $\mathcal N$ taking values exclusively in ${\mathcal M} \subseteq
{\mathcal M}'$. This showes the existence of the adjoint operator $T^*$
of $T$ in the sense of its definition.

The property of the graph of a bounded module operator to be a topological
summand with topological complement $\{ (0,z) : z \in \mathcal N \}$ follows
from the decomposition $(x,y) = (x,T(x)) + (0,y-T(x))$ for every $x \in
\mathcal M$, $y \in \mathcal N$. Since $T(0)=0$ for any linear operator $T$
the intersection of the graph with the $A$-$B$ submodule $\{ (0,z) : z \in
\mathcal N \}$ is always trivial.
\end{proof}

\begin{problem}
  It would be highly interesting to know whether the kernel
  of every surjective bounded module operator would be merely a topological
  summand, or whether there are counterexamples.
\end{problem}

The observation above gives us the opportunity to add two further equivalent
conditions to Theorem \ref{Mag-Schw}.
Much less obvious is the fact that only topological summands have to be
considered to characterize the same class of coefficient C*-algebras as a
distinguished one in the research field of Hilbert C*-modules. This is
established by the next theorem which also demonstrates the missing in
\cite{FrTr} self-duality assumption to be inevitable for obtaining the results
of \cite[\S 1]{FrTr} (cf.~\cite{Manuilov:00}).

\begin{theorem}  \label{Fr-neu2}
   There are further equivalent conditions to the conditions listed in Theorem
   \ref{Mag-Schw} and Corollary \ref{Fr-neu}:

   \begin{list}{(\roman{cou001})}{\usecounter{cou001}}
   \setcounter{cou001}{5}
   \item  The kernels of all bounded $A$-linear operators between arbitrary
          Hilbert $A$-mo\-du\-les are orthogonal summands.
   \item  The images of all bounded $A$-linear operators with norm-closed
          range between arbitrary Hilbert $A$-modules are orthogonal summands.
   \item  For every Hilbert $A$-module every Hilbert $A$-submodule
          is automatically topologically complemented there, i.e.~it is a
          topological summand.
   \item  For every (maximal) norm-closed left ideal $I$ of $A$ the
          corresponding open projection $p \in A^{**}$ is an element of the
          multiplier C*-algebra ${\rm M}(A)$ of $A$.
   \end{list}
\end{theorem}

\begin{proof}
By Theorem \ref{Mag-Schw}, (ii) condition (i) on $A$ implies condition (viii).
Since norm-closed one-sided ideals $I$ of $A$ and open projections
\linebreak[4]
$p \in A^{**}$ are in one-to-one correspondence by \cite{Ped} a close
examination of the norm-closed left ideals of C*-algebras $A$ with property
(i) shows assertion (ix) to be valid.

If condition (ix) on $A$ holds then the multiplier C*-algebra
${\rm M}(A)$ is an atomic type I von Neumann algebra, cf.
\cite[Lemma 2]{Lin1}. It can be represented as a direct sum over a
discrete measure space of a certain cardinality, say
$l_\infty$-$\sum_i {\rm B}(H_i)$ for Hilbert spaces $H_i$ by the
direct integral decomposition theory of von Neumann algebras. If
one of the Hilbert spaces $H_i$ is infinite-dimensional the
appropriate $i$-th block in the direct sum decomposition of $A$
above cannot coincide with the entire set ${\rm B}(H_i)$ since the
latter contains maximal one-sided ideals not supported by
projections of ${\rm B}(H_i)$ itself. So the $i$-th block has to
be $*$-isomorphic to ${\rm K}(H_i)$ for every Hilbert space $H_i$.
Considering the two-sided ideal $\overline{c_0A}$ in $A$ the
corresponding carrier projection is central in $A$, and so we can
resort to the center of $A$ to continue our considerations. The
center of $A$ cannot be equal to $l_\infty$ since the norm-closed
ideal $c_0$ of $l_\infty$ would not be supported by a projection
inside $l_\infty$ as supposed. So the center of $A$ has to be
either finite-dimensional or a $c_0$-space of arbitrary
cardinality. This implies condition (i).

Furthermore, suppose the C*-algebra $A$ fulfills condition (viii)
and $A \not\equiv {\mathbb C}$. Consider $A$ and an arbitrary
maximal left ideal $I$ of $A$ as Hilbert $A$-modules equipped with
the standard $A$-valued inner product of the C*-algebra $A$. Since
$I$ is a Hilbert $A$-submodule of $A$ and all bounded module maps
on $A$ can be identified with right multipliers ${\rm RM}(A)$ of
$A$ by \cite{Fr98} there exists an idempotent right multiplier $p
\not= 1_A$ of $A$ such that $x=xp$ for any $x \in I$ by
assumption. Indeed, the element $t = pp^* + (1-p)^*(1-p)$ is
invertible in $M(A)$, and $q = t^{-1}pp^* =pp^*t^{-1} \in M(A)$ is
the orthogonal projection with the same range as $p$. Note, that
$p$ and $q$ commute. Since the orthogonal complement of the
(positive) carrier projection $q \in A^{**}$ of $I$ is a
non-trivial minimal projection of $A^{**}$ and $p$ and $q$ commute
in $A^{**}$, the coincidence $p=q$ and the inclusion $q \in {\rm
M}(A)$ turn out to hold. We arrive at assertion (ix) because the
maximal ideal $I$ of $A$ was assumed to be an arbitrary one and
any other norm-closed left ideal equals to the intersection of
maximal left ideals. The established equivalence of (ix) and (i)
completes the argument.

Consider the situation of a graph of a bounded $A$-linear operator
$T: \mathcal M \to \mathcal N$ between two Hilbert $A$-modules. By Theorem
\ref{Mag-Schw}, Corollary \ref{Fr-neu} and Corollary \ref{adjointability}
the operator $T$ is adjointable if and only if its graph is an orthogonal
summand. However, in any case its graph is the kernel of the bounded
$A$-linear operator $S: (x,y) \to (0,T(x)-y)$ that acts on $\mathcal M \oplus
\mathcal N$. Consequently, condition (vi) is equivalent to the adjointability
of all bounded module operators between Hilbert $A$-modules. Resorting to the
special case of canonical embeddings of maximal ideals $I$ of $A$ into $A$
we obtain the equivalence of (vi) with (ix) and hence, with Theorem
\ref{Mag-Schw}, (i).

If a bounded $A$-linear operator $T: \mathcal M \to \mathcal N$
between two Hilbert $A$-modules admits an image which is an
orthogonal summand then the graph of $T$ is the image of the
bounded module operator $({\rm id}_{\mathcal M} \oplus
(0_{\mathcal N} + T))$ on  $\mathcal M \oplus \mathcal N$.
Moreover this set is an orthogonal summand of $\mathcal M \oplus
\mathcal N$, and it is the graph of $T$. Again, the graph of $T$
is an orthogonal summand of $\mathcal M \oplus \mathcal N$ if and
only if $T$ is adjointable, and the arguments above can be
repeated demonstrating the equivalence of (vii) and Theorem
\ref{Mag-Schw}, (i).
\end{proof}

\begin{problem} \label{problemXX}
 Is any Hilbert C*-submodule of a any Hilbert C*-module $\mathcal H$
 that coincides with its bi-orthogonal complement inside $\mathcal H$, the
 kernel of a bounded $A$-linear operator mapping $\mathcal H$ into itself (or,
 alternatively, to other Hilbert C*-modules)?
\end{problem}

\begin{problem} \label{problemYY}
 Characterize those C*-algebras $A$ for which the following condition holds:
  For every Hilbert C*-module over $A$ every Hilbert $A$-submodule which
  coincides with its bi-ortho\-gonal complement is automatically
  topologically complemented there.
\end{problem}

Problem \ref{problemYY} remodels the difference between
B.~Magajna's theorem and J.~Schwei\-zer's theorem on the level of
topological summands. The results of M. Kusuda \cite{Kusuda2}
indicate that the solution has to be similar to that one given on
the level of orthogonal summands. However, research has to be
continued.

\bigskip \noindent
{\bf Acknowledgements:} I am indebted to V.~I.~Paulsen and
D.~P.~Blecher for their hospitality and support during my one-year
stay at the University of Houston in Houston, Texas, U.S.A., in
1998. The fruitful discussions on injectivity and projectivity of
Hilbert C*-modules and C*-algebras Vern I.~Paulsen and I had that
time stimulated also the research the results of which are
contained in the present paper. Also, I am grateful to the referee
for his substantial remarks which improved the paper in places
substantially.



\end{document}